\documentclass[11pt,centertags,twoside]{amsart}
\usepackage{amssymb}
\usepackage[mathcal]{euscript}
\usepackage{mathptm}

\DeclareSymbolFont{SY}{U}{psy}{m}{n}
\DeclareMathSymbol{\emptyset}{\mathord}{SY}{'306}

\renewcommand{\epsilon}{\varepsilon}

\numberwithin{equation}{section}

\newcommand{\C}{\mathbb{C}}

\newcommand{\N}{\mathbb{N}}

\newcommand{\R}{\mathbb{R}}
\newcommand{\EE}{\mathsf{E}}

\newcommand{\bA}{{\mathbf A}}
\newcommand{\bB}{{\mathbf B}}

\newcommand{\bV}{{\mathbf V}}

\newcommand{\cB}{{\mathcal B}}

\newcommand{\cG}{{\mathcal G}}

\newcommand{\cL}{{\mathcal L}}

\newcommand{\cP}{{\mathcal P}}

\newcommand{\cS}{{\mathcal S}}
\newcommand{\cT}{{\mathcal T}}

\renewcommand{\Re}{{\ensuremath{\mathrm{Re}}}}

\newcommand{\dist}{{\ensuremath{\mathrm{dist}}}}

\newcommand{\fH}{\mathfrak{H}}

\newcommand{\fK}{\mathfrak{K}}
\newcommand{\fL}{\mathfrak{L}}
\newcommand{\fM}{\mathfrak{M}}
\newcommand{\fN}{\mathfrak{N}}

\newcommand{\fQ}{\mathfrak{Q}}

\DeclareMathOperator*{\slim}{s-lim}

\DeclareMathOperator{\Ran}{\mathrm{Ran}}
\DeclareMathOperator{\Ker}{\mathrm{Ker}\,}

\DeclareMathOperator{\Dom}{\mathrm{Dom}}
\newcommand{\spec}{{\ensuremath{\mathrm{spec}}}}

\newcommand{\specpp}{\spec_{\mathrm{pp}}}

\newtheorem{theorem}{Theorem}[section]
\newtheorem{lemma}[theorem]{Lemma}

\newtheorem{corollary}[theorem]{Corollary}
\newtheorem{hypothesis}[theorem]{Hypothesis}

\theoremstyle{remark}
\newtheorem{remark}[theorem]{Remark}{\bf}{\rm}
{\bf}{\rm}

\title[The $\tan 2\Theta$ Theorem]
{A generalization of the tan\,2$\Theta$ Theorem}

\author[V.~Kostrykin]{Vadim Kostrykin}
\address[Vadim Kostrykin]{Fraunhofer-Institut f\"{u}r
Lasertechnik, Steinbachstra{\ss}e 15, D-52074\\ Aachen, Germany}
\email{kostrykin@ilt.fraunhofer.de, kostrykin@t-online.de}

\author[K.~A.~Makarov]{Konstantin A.~Makarov}
\address[Konstantin A. Makarov]{Department of Mathematics, University of
Missouri, Co\-lum\-bia, MO 65211, USA}
\email{makarov@math.missouri.edu}

\author[A.~K.~Motovilov]{Alexander K.~Motovilov}
\address[Alexander K.~Motovilov]{Department of Mathematics, University of
Missouri, Co\-lum\-bia, MO 65211, USA\newline \textit{Permanent Address:
Joint Institute for Nuclear Research, 141980 Dubna, Moscow Region, Russia} }
\email{motovilv@thsun1.jinr.ru}



\date{February 1, 2003}

\begin{document}

\begin{abstract}
Let $\bA$ be a bounded self-adjoint operator on a separable Hilbert space
$\fH$ and  $\fH_0\subset\fH$  a closed invariant subspace of $\bA$.
Assuming that $\sup\spec(A_0)\leq \inf\spec(A_1)$, where $A_0$ and $A_1$
are restrictions of $\bA$ onto the subspaces $\fH_0$ and
$\fH_1=\fH_0^\perp$, respectively, we study the variation of the invariant
subspace $\fH_0$ under bounded self-adjoint perturbations $\bV$  that are
off-diagonal with respect to the decomposition $\fH=\fH_0\oplus\fH_1$. We
obtain sharp two-sided estimates on the norm of the difference of the
orthogonal projections onto invariant subspaces of the operators $\bA$ and
$\bB=\bA+\bV$. These results extend the celebrated Davis-Kahan $\tan
2\Theta$ Theorem. On this basis we also prove new existence and uniqueness
theorems for contractive solutions to the operator Riccati equation, thus,
extending recent results of Adamyan, Langer, and Tretter.
\end{abstract}

\maketitle

\section{Introduction}\label{sec:1}

Given a self-adjoint bounded operator $\bA$ and a closed invariant subspace
$\fH_0\subset\fH$ of $\bA$ we set $A_i=\bA|_{\fH_i}$, $i=0,1$ with
$\fH_1=\fH\ominus\fH_0$. Assuming that the perturbation $\bV$ is
off-diagonal with respect to the orthogonal decomposition $\fH=\fH_0 \oplus
\fH_1$ consider the $2\times 2$ self-adjoint  operator matrix
\begin{equation*}
\bB=\bA+\bV=\begin{pmatrix} A_0 & V \\ V^\ast & A_1
\end{pmatrix},
\end{equation*}
where $V$ is a bounded operator from $\fH_1$ to $\fH_0$.

In the 1970 paper \cite{Davis:Kahan} Davis and Kahan proved
that if
\begin{equation}\label{eka}
\sup\spec(A_0) < \inf  \spec(A_1),
\end{equation}
then the difference of the spectral projections
\begin{equation*}
P=\EE_{\bA}\big((-\infty,\sup\spec(A_0)]\big)\quad\text{and}\quad
Q=\EE_{\bB}\big((-\infty,\sup\spec(A_0)]\big)
\end{equation*}
for the operators $\bA$ and $\bB$, respectively, corresponding
to the interval $(-\infty,\sup\spec(A_0)]$ admits the estimate
\begin{equation}\label{tan:2:Theta}
  \|P-Q\|\leq \sin \bigg (
\frac{1}{2}\arctan \frac{2\|V\|}{d} \bigg
)<\frac{\sqrt{2}}{2},
\end{equation}
where
\begin{equation*}
d=\dist (\spec (A_0) , \spec(A_1)).
\end{equation*}
Estimate \eqref{tan:2:Theta} can be equivalently expressed as
the $\tan 2\Theta$ Theorem:
\begin{equation*}
\|\tan 2\Theta\| \leq \frac{2\|V\|}{d},\qquad
\spec(\Theta)\subset [0,\pi/4),
\end{equation*}
where $\Theta$ is the operator angle between the subspaces
$\Ran P$ and $\Ran Q$ (see, e.g.,
\cite{Kostrykin:Makarov:Motovilov:2}).

By known results on graph subspaces (see, e.g.,
\cite{Albeverio}, \cite{Apostol:Foias:Salinas},
\cite{Daughtry}, \cite{Kostrykin:Makarov:Motovilov:2})
estimate \eqref{tan:2:Theta} in particular implies that  the
Riccati equation
\begin{equation}\label{intro:Riccati}
A_1 X - X A_0 - XVX + V^\ast = 0
\end{equation}
has a contractive solution $X:\fH_0\rightarrow\fH_1$
satisfying the norm estimate
\begin{equation}\label{tan:2Theta}
\|X\| = \frac{\|P-Q\|}{\sqrt{1-\|P-Q\|^2}}\leq \tan
\bigg(\frac{1}{2}\arctan \frac{2\|V\|}{d} \bigg )<1.
\end{equation}
Moreover, the graph of $X$, i.e., the subspace $\cG(\fH_0,X):=\{x\oplus Xx|
x\in\fH_0\}$, coincides with the spectral subspace $\Ran
\EE_{\bB}\big((-\infty, \sup\spec(A_0)]\big)$ of the operator $\bB$.

Independently of the work of Davis and Kahan the existence of a unique
contractive solution to the Riccati equation  under condition \eqref{eka}
has been proven by Adamyan and Langer in \cite{Adamyan:Langer:95}, where
the operators $A_0$ and $A_1$ were allowed to be semibounded. In a recent
paper by Adamyan, Langer, and Tretter \cite{Adamyan:Langer:Tretter:2000a}
the existence result has been extended to the case where the spectra of
$A_0$ and $A_1$ intersect at one point $\lambda\in \R$, that is,
\begin{equation*}
\sup \spec(A_0) =\inf \spec(A_1)=\lambda,
\end{equation*}
provided that at least one of the following conditions
\begin{equation}\label{1}
\Ker (A_0-\lambda)=\{0\}, \quad \Ker (A_1-\lambda)\cap \Ker
V=\{0\}
\end{equation}
or
\begin{equation}\label{2}
\Ker (A_1-\lambda)=\{0\}, \quad \Ker (A_0-\lambda)\cap \Ker
V^\ast =\{0\}
\end{equation}
holds. In this case the Riccati equation \eqref{intro:Riccati} has been
proven to have a unique contractive solution $X$, which appears to be a
strict contraction. The graph of $X$, as above, coincides with  the
spectral subspace $\Ran \EE_\bB((-\infty, \lambda))$ of the operator $\bB$.

Conditions \eqref{1} and \eqref{2} are rather restrictive. In
particular, in this case $\lambda$ may  be an eigenvalue
neither for both $A_0$ and $A_1$ nor for $\bB$.

The main goal of the present article is to drop  conditions
\eqref{1} and \eqref{2} and to carry out the analysis under
the only assumption that
\begin{equation}\label{only}
\sup\spec(A_0)\leq \lambda\leq \inf\spec(A_1).
\end{equation}

Below we will prove (see Theorem \ref{Davis:Kahan}) that under
hypothesis \eqref{only} the $\bB$-invariant subspace
\begin{equation}\label{intro:Q}
\fQ=\Ran \EE_\bB((-\infty, \lambda))\oplus\left(\Ker
(A_0-\lambda)\cap \Ker V^\ast\right)
\end{equation}
is the graph of a contractive operator $X:\fH_0\rightarrow\fH_1$. Moreover,
the norm of the operator $X$ satisfies the lower bound
\begin{equation*}
\|X\| \geq \frac{\delta}{\|V\|},
\end{equation*}
where
\begin{equation*}
\delta = \max\{\inf\spec(\bA) - \inf\spec(\bB), \sup\spec(\bB) - \sup\spec(\bA)\} \geq 0
\end{equation*}
is the maximal shift of the edges of the spectrum of the
operator $\bA$ under the perturbation $\bV$. These results can
be stated equivalently as the two-sided estimate
\begin{equation}\label{first:inequality}
\frac{\delta}{\sqrt{\delta^2 +\|V\|^2}} \leq \| P-Q\| \leq \frac{\sqrt{2}}{2},
\end{equation}
where $P$ and $Q$ are orthogonal projections in $\fH$ onto the subspaces
$\fH_0$ and $\fQ$, respectively. Notice that the for the subspace $\fQ$ to
be a spectral subspace of $\bB$ it is necessary and sufficient
 that
\begin{equation*}
\text{either}\quad\Ker(A_0-\lambda)\cap\Ker
V^\ast=\{0\}\quad\text{or}\quad \Ker(A_1-\lambda)\cap\Ker V
=\{0\},
\end{equation*}
and then, necessarily,
\begin{equation*}
\text{either}\quad \fQ=\Ran \EE_\bB((-\infty, \lambda)) \quad \text{or}
\quad \fQ=\Ran \EE_\bB((-\infty, \lambda]).
\end{equation*}

The fact that the subspace $\fQ$ is a graph of a contractive operator $X$
means that the Ric\-ca\-ti equation \eqref{intro:Riccati} has a contractive
solution. In contrast to the case studied in
\cite{Adamyan:Langer:Tretter:2000a}, the solution $X$ is in general (under
hypothesis \eqref{only}) neither necessarily strictly contractive nor
unique in the set of all contractive solutions to the Riccati equation.
Moreover, if the invariant subspace $\fQ$ is not a spectral subspace of the
operator $\bB$, then $X$ is a non-isolated point (in the operator norm
topology) of the set of all solutions of \eqref{intro:Riccati}.

Below we will prove (see Theorem \ref{Clubshilfe}) that under
assumption \eqref{only} the operator $X$ is the unique
solution to the Riccati equation \eqref{intro:Riccati} within
the class of bounded linear operators from $\fH_0$ to $\fH_1$
satisfying the additional requirements that
\begin{equation*}
\Ker(A_0-\lambda)\cap \Ker V^\ast \subset \Ker
X\quad\text{and}\quad \spec(A_0 + V
X)\subset(-\infty,\lambda].
\end{equation*}
{}Furthermore, we formulate and prove necessary and sufficient
conditions for a contractive solution to be  a  unique
contractive solution (see Theorem \ref{4.11a} or/and  Theorem
\ref{4.11}). The solution $X$ satisfying $\cG(\fH_0,X)=\fQ$ is
shown to be the unique contractive solution to
\eqref{intro:Riccati} if and only if it is strictly
contractive (see Corollary \ref{cor:strictly:contractive}).
Note that in the case where  the contractive solution is
non-unique but its graph is a spectral subspace for the
operator $\bB$, a complete description of the set of all
contractive solutions to the Riccati equation can be given by
means of Theorem 6.2 in \cite{Kostrykin:Makarov:Motovilov:2}.

The technics developed in the present work to prove that the subspace $\fQ$
\eqref{intro:Q} is the graph of a contractive operator (Theorem
\ref{Davis:Kahan}), is an extension of the geometric ideas of Davis and
Kahan \cite{Davis:123}, \cite{Davis:Kahan}. The existence and uniqueness
results (Theorems \ref{Clubshilfe}, \ref{4.11a}, and \ref{4.11}) are
obtained in the framework of the geometric approach of our recent paper
\cite{Kostrykin:Makarov:Motovilov:2}. The previously known results by Davis
and Kahan \cite{Davis:123}, \cite{Davis:Kahan}, Adamyan and Langer
\cite{Adamyan:Langer:95}, and Adamyan, Langer, and Tretter
\cite{Adamyan:Langer:Tretter:2000a} appear to be their direct corollaries.

A few words about the notations used throughout the paper. Given a linear
operator $A$ on a Hilbert space $\fK$, by $\spec(A)$ we denote the spectrum
of $A$. If not explicitly stated otherwise, $\fN^\perp$ denotes the
orthogonal complement in $\fK$ of a subspace $\fN\subset\fK$, i.e.,
$\fN^\perp=\fK\ominus\fN$. The identity operator on $\fK$ is denoted by
$I_\fK$. The notation $\cB(\fK,\fL)$ is used for the set of bounded
operators from the Hilbert space $\fK$ to the Hilbert space $\fL$. Finally,
we write $\cB(\fK)=\cB(\fK,\fK)$.

\subsection*{Acknowledgments.} V.~Kostrykin is grateful to V.~Enss,
A.~Knauf, and R.~Schrader for useful discussions. A.~K.~Motovilov
acknowledges the great hospitality and financial support by the Department
of Mathematics, University of Missouri--Columbia, MO, USA. He was also
supported in part by the Russian Foundation for Basic Research within
Project RFBR 01-01-00958.

\section{Upper Bound}\label{S2island}
\setcounter{equation}{0}

Throughout the whole work we adopt the following hypothesis.

\begin{hypothesis}\label{diag}
Assume that the separable Hilbert space $\fH$ is decomposed
into the orthogonal sum of two subspaces
\begin{equation}\label{decom}
\fH=\fH_0 \oplus \fH_1.
\end{equation}
Assume, in addition, that $\bB$ is a self-adjoint operator on $\fH$
 represented with respect to the decomposition
\eqref{decom} as a $2\times2$ operator block matrix
\begin{equation*}
\bB=\begin{pmatrix}
  A_0 & V \\
  V^\ast & A_1
\end{pmatrix},
\end{equation*}
where $A_i$ are bounded self-adjoint operators in $\fH_i$,
$i=0,1$, while $V$ is a bounded operator from $\fH_1$ to
$\fH_0$. More explicitly, $\bB =\bA+\bV$, where $\bA$ is the
bounded diagonal self-adjoint operator,
\begin{equation*}
\bA = \begin{pmatrix}
  A_0 & 0 \\
  0 & A_1
\end{pmatrix},
\end{equation*}
and the operator $\bV=\bV^\ast$ is an off-diagonal bounded
operator
\begin{equation*}
\bV=\begin{pmatrix}
  0 & V \\
  V^\ast & 0
\end{pmatrix}.
\end{equation*}
Moreover, assume that
\begin{equation}\label{only:2}
\sup\spec(A_0)\leq \lambda\leq \inf\spec(A_1)
\end{equation}
for some $\lambda\in \R$\,.
\end{hypothesis}

If, under Hypothesis \ref{diag},  $\lambda$ is a multiple
eigenvalue of the operator $\bB$, then $\bB$ has infinitely
many invariant subspaces $\fL_{\mathbf{B}}$  such that
\begin{equation}\label{88}
\Ran\EE_{\bB}((-\infty, \lambda)) \subsetneq
 \fL_{\bB}\subsetneq
\Ran\EE_{\bB}((-\infty, \lambda])
\end{equation}
that are necessarily not spectral subspaces.

A  criterion for a $\bB$-invariant subspace $\fL_{\bB}$
satisfying \eqref{88} to be a graph subspace associated with the
decomposition $\fH=\fH_0\oplus \fH_1$ is given by the following theorem.

\begin{theorem}\label{n0wings}
Assume Hypothesis \ref{diag}. Then
\begin{equation}
\label{canon} \Ker(\bB-\lambda) = (\Ker(A_0-\lambda)\cap \Ker
V^\ast)\oplus (\Ker(A_1-\lambda)\cap \Ker V).
\end{equation}
Moreover, given a subspace $\fN\subset \Ker(\bB-\lambda)$, the
subspace
\begin{equation}\label{graf11}
\fL_{\bB}:=\Ran\EE_{\bB}((-\infty, \lambda))\oplus \fN
\end{equation}
is a graph subspace associated with the subspace $\fH_0$ in the
decomposition $\fH=\fH_0\oplus \fH_1$ if and only if $\fN$ is a
graph subspace associated with the subspace $\Ker(A_0-\lambda)\cap
\Ker V^\ast$ in the
decomposition \eqref{canon}.
\end{theorem}

\begin{proof}
Without loss of generality we may set $\lambda=0$. The
inclusion
\begin{equation}\label{obratno}
(\Ker A_0 \cap \Ker V^\ast)\oplus (\Ker A_1 \cap \Ker
V)\subset \Ker \bB
\end{equation}
is obvious. In order to prove the opposite inclusion assume
that $\begin{pmatrix} x\\ y
\end{pmatrix}\in\Ker \bB$ with $x\in\fH_0$ and $y\in\fH_1$, i.e.,
\begin{equation}\label{ref}
A_0 x + V y = 0\qquad \text{ and } \qquad  V^\ast x + A_1 y =
0.
\end{equation}
Suppose that $x\notin\Ker A_0$. Then $(x,A_0 x)<0$ and
therefore $(x,Vy)>0$ using the first equation in \eqref{ref}.
{}From the second of equations \eqref{ref} it follows that
$(y, A_1y)<0$, which is in a contradiction with
\eqref{only:2}.
 Thus $x\in\Ker A_0$. Similarly one proves that $y\in\Ker
A_1 $. By using \eqref{ref} it follows that $x\in\Ker V^\ast$
and $y\in\Ker V$  which together with \eqref{obratno} proves
\eqref{canon}.

In order to prove the second statement of the theorem notice
that by Theorem \ref{3.2} in the Appendix the subspace
$\fN\subset \Ker \bB$ is a graph subspace associated with the
subspace $\Ker A_0\cap\Ker V^\ast$ in the decomposition
\eqref{canon} (recall that we assumed that $\lambda=0$) if and
only if
\begin{equation*}
(\Ker A_0\cap \Ker V^\ast) \cap (\Ker \bB\ominus \fN)= (\Ker
A_1 \cap \Ker V) \cap  \fN=\{0\}.
\end{equation*}
Again from Theorem \ref{3.2}  it follows that the subspace
$\Ran\EE_{\bB}((-\infty,0))\oplus \fN$ is a graph subspace
associated with the subspace $\fH_0$ in the decomposition
$\fH=\fH_0\oplus \fH_1$ if and only if
\begin{equation}
\fH_0\cap \big(\Ran \EE_{\bB}((-\infty, 0))\oplus
\fN\big)^\perp =
 \fH_1\cap
(\Ran \EE_{\bB}((-\infty, 0))\oplus \fN)=\{0\}.
\end{equation}
Therefore, to complete the proof it is sufficient to establish
the following equalities
\begin{equation}\label{mnmn}
\fH_0\cap \big(\Ran \EE_{\bB}((-\infty, 0))\oplus
\fN\big)^\perp = (\Ker A_0\cap \Ker V^\ast) \cap (\Ker
\bB\ominus \fN)
\end{equation}
and
\begin{equation}\label{nmnm}
\fH_1\cap \big(\Ran \EE_{\bB}((-\infty, 0))\oplus
\fN\big)=(\Ker A_1\cap \Ker V)\cap\fN.
\end{equation}

First, we prove that the left-hand side of \eqref{mnmn} is a
subset of the right-hand  side of
\eqref{mnmn}, i.e.,
\begin{equation}\label{inin}
\fH_0\cap \big(\Ran \EE_{\bB}((-\infty, 0))\oplus
\fN\big)^\perp \subset (\Ker A_0\cap \Ker V^\ast) \cap (\Ker
\bB\ominus \fN).
\end{equation}
Let $x\in \fH_0 $ and $x\perp \Ran \EE_{\bB}((-\infty,
0))\oplus \fN$. Clearly,
\begin{equation*}
(x, \bB x)\geq 0\quad \text{ for }\quad x\perp \Ran
\EE_{\bB}((-\infty, 0 ))\oplus \fN.
\end{equation*}
Moreover,
\begin{equation*}
(x, \bB x)=(x, \bA x)\leq 0 \quad \text{ for all }\quad x\in
\fH_0,
\end{equation*}
since the operator matrix $\mathbf{V}$ is off-diagonal and the
subspace $\fH_0$ is $\mathbf{A}$-invariant. Hence,
\begin{equation*}
(x, \bA x)=(x, \bB x)=0\quad \text{ for } x\in \fH_0 \cap
\big(\Ran \EE_{\bB}((-\infty,0))\oplus \fN \big)^\perp,
\end{equation*}
which implies $x\in \Ker \bB\ominus \fN$ by the variational principle.
Therefore,  the inclusion \eqref{inin} holds.

The opposite inclusion
\begin{equation*}
\fH_0\cap \big(\Ran \EE_{\bB}((-\infty, 0))\oplus
\fN\big)^\perp \supset (\Ker A_0\cap \Ker V^\ast) \cap (\Ker
\bB\ominus \fN)
\end{equation*}
is obvious, which proves \eqref{mnmn}.

The equality \eqref{nmnm} is proven in a similar way.
\end{proof}
\begin{remark}\label{remark:on:Q}
As it follows from \eqref{canon} the closed $\bB$-invariant subspace
\begin{equation}\label{Q:Def}
\fQ=\Ran\EE_{\bB}((-\infty,\lambda))\oplus
\big(\Ran\EE_{\bB}(\{\lambda\})\cap \fH_0\big)\subset\fH.
\end{equation}
 is a spectral
subspace for the operator $\bB$ if and only if
\begin{equation*}
\text{either}\quad \Ker(A_0-\lambda)\cap\Ker V^\ast =
\{0\}\quad \text{or}\quad\Ker(A_1-\lambda)\cap\Ker V = \{0\}.
\end{equation*}
\end{remark}

The following theorem  characterizes the subspace $\fQ$ as the
graph  of some contractive operator $X$ from $\fH_0$
to $\fH_1$, that is, $\fQ=\cG(\fH_0,X)$ where
\begin{equation*}
\cG(\fH_0,X) = \{x\oplus Xx| x\in\fH_0\}.
\end{equation*}

\begin{theorem}\label{Davis:Kahan}
Assume Hypothesis \ref{diag}. Then:

\emph{(i)} The  $\bB$-invariant  subspace $\fQ$
 is a graph subspace $\cG(\fH_0, X)$
associated with the subspace $\fH_0$ in decomposition
\eqref{decom}, where the operator $X$ is a contraction,
$\|X\|\leq 1$, with the properties
\begin{equation}\label{leer}
\Ker(A_0-\lambda) \cap \Ker V^\ast  \subset \Ker X, \qquad
\Ker(A_1-\lambda) \cap \Ker V  \subset \Ker X^\ast.
\end{equation}

\emph{(ii)} If either $\Ker (A_0-\lambda)=\{0\}$ or $\Ker
(A_1-\lambda)=\{0\}$, then $X$ is a strict contraction, i.e.,
\begin{equation*}
\|Xf\| < \|f\|, \quad f\ne 0.
\end{equation*}

\emph{(iii)} If
\begin{equation}\label{razdelenu}
d=\dist(\spec(A_0), \spec (A_1))>0,
\end{equation}
then  $X$ is a uniform contraction satisfying the estimate
\begin{equation}\label{normasnizu}
\|X\|\le \tan \left(\frac{1}{2}\arctan \frac{2\|V\|}{d}\right)
< 1.
\end{equation}
\end{theorem}

\begin{remark}
In Section \ref{sec:unique} we will establish necessary and sufficient
conditions guaranteeing that the operator $X$ referred to in Theorem
\ref{Davis:Kahan} is a strict contraction (see Corollary
\ref{cor:strictly:contractive} below). The claim (ii) of Theorem
\ref{Davis:Kahan} will then appear to be a corollary of this more general
result.

\end{remark}

\begin{proof}[Proof of Theorem \ref{Davis:Kahan}]
Without loss of generality we assume that $\lambda=0$.

(i) Step 1.  First we prove the assertion under the additional
assumption
 that $\Ker \bB=\{0\}$.

By Theorem \ref{n0wings} the subspace  $\fQ$ is a graph subspace associated
with the decomposition $\fH=\fH_0\oplus\fH_1$, i.e.,
\begin{equation}\label{graf}
\fQ =\cG(\fH_0, X),
\end{equation}
where $X$ is a (possibly unbounded) densely defined closed
operator from $\fH_0$ to $\fH_1$. Let $X=U|X|$ be the polar
decomposition for $X$, where $U:\fH_0\to \fH_1$ is a partial
isometry with the initial subspace $(\Ker X)^\perp$ and the
final subspace $\overline{\Ran X}$ and $|X|=(X^\ast X)^{1/2}$,
the absolute value of $X$.

First we show that
\begin{equation}\label{net}
\specpp(|X|)\subset [0,1],
\end{equation}
where  $\specpp(|X|)$ denotes the set of all eigenvalues of
$|X|$.

Let $0\ne \mu\in \specpp(|X|)$  and $f$ be an eigenvector of
$|X|$ corresponding to the eigenvalue $\mu$, i.e.,
\begin{equation}\label{sobstv}
|X|f=\mu f, \quad 0\ne f\in \Dom (|X|).
\end{equation}
By \eqref{graf}
\begin{equation}\label{F}
F= f\oplus X f=f\oplus \mu \, U f\in \fQ
\end{equation}
using \eqref{sobstv}. Since $f \perp \Ker |X|=\Ker X$, the
element $f$ belongs to the initial subspace of the isometry
$U$. Moreover,
\begin{equation*}
U f=\mu^{-1}U|X|f\in \Ran X \subset \overline{\Ran X},
\end{equation*}
i.e., the element $U f$ belongs to the final subspace of $U$
and hence $U^\ast Uf=f$, which, in particular, proves that
$Uf\in \Dom (X^\ast)$. Therefore,
\begin{equation}\label{G}
G= (-X^\ast Uf)\oplus Uf=(-\mu\, f )\oplus Uf \in \fQ^\perp.
\end{equation}

Using \eqref{F}, \eqref{G}, and the hypothesis $\Ker \bB
=\{0\}$ one obtains the following two strict inequalities
\begin{alignat}{3}
0\ & >\ (F, \bB F)\ & =\ &\ (f, A_0 f) + 2\mu\, \Re (V^\ast f,
Uf) + \mu^2(Uf, A_1Uf),
\label{bole}\\
0\ & <\  (G, \bB G)\ & =\ &\ \mu^2(f, A_0 f)-2\mu\, \Re
(V^\ast f, Uf) +(Uf, A_1 Uf). \label{menee}
\end{alignat}

If $\mu>0$ satisfies \eqref{bole} and \eqref{menee}, then necessarily
$\mu\le 1$. In order to see that we subtract \eqref{menee} from
\eqref{bole} getting the inequality
\begin{equation}\label{ner}
(1-\mu^2)\big((Uf, A_1 Uf)-(f,A_0 f)\big) > 4\mu \, \Re
(V^\ast f, Uf).
\end{equation}
Since $A_0\leq 0$ and $A_1\geq 0$, equation \eqref{ner}
implies $\Re (V^\ast f, Uf)<0$ for $\mu> 1$ which contradicts
the orthogonality of the elements $F$ and $\bB G$:
\begin{equation}\label{ravno}
 (F, \bB G)= \mu\big((Uf,A_1 Uf)-(f,A_0 f)\big)+(1-\mu^2)\Re
(V^\ast f, Uf)=0.
\end{equation}
Hence, \eqref{net} is proven.

Our next goal is to prove that the operator $X$ is a
contraction.

Let $\{\cP_n^{(0)}\}_{n\in\N}$ and $\{\cP_n^{(1)}\}_{n\in\N}$
be two sequences of finite-dimensional orthogonal projections
such that $\Ran \cP_n^{(0)}\subset \fH_0$,
$\Ran\cP_n^{(1)}\subset \fH_1$, and
\begin{equation}\label{plimits}
\slim_{n\rightarrow\infty}\cP_n^{(0)} = P, \qquad
\slim_{n\rightarrow\infty}\cP_n^{(1)} = P^\perp,
\end{equation}
where $P$ is the orthogonal projection  from $\fH$ onto
$\fH_0$ and
\begin{equation}\label{limits:0}
\slim_{n\rightarrow\infty} \EE_{\mathbf{A}_n}(\{0\}) =
\EE_{\mathbf{A}}(\{0\}),
\end{equation}
\begin{equation}\label{limits:1}
\slim_{n\rightarrow\infty} \EE_{\mathbf{A}_n}((-\infty,0)) =
\EE_{\mathbf{A}}((-\infty, 0)),
\end{equation}
where
\begin{equation*}
\mathbf{A}_n = \begin{pmatrix}\cP_n^{(0)} A_0 \cP_n^{(0)} & 0
\\ 0 & \cP_n^{(1)}A_1\cP_n^{(1)} \end{pmatrix}
\end{equation*}
are the corresponding finite-dimensional truncations of the operator $\bA$.
The existence of such sequences can easily be shown by
splitting off the subspaces $\Ker A_0$ and $\Ker A_1$.

Introducing the finite rank operators
\begin{equation*}
\mathbf{V}_n=\begin{pmatrix}0 &
\cP_n^{(0)}V\cP_n^{(1)} \\
\cP_n^{(1)} V^\ast \cP_n^{(0)} & 0 \end{pmatrix}
\end{equation*}
one concludes (see, e.g., Theorem I.5.2 in
\cite{Birman:Solomyak})
 that
\begin{equation}\label{slim}
\slim_{n\rightarrow\infty}(\mathbf{A}_n+\mathbf{V}_n) = \bB.
\end{equation}
Since $\Ker \bB=\{0\}$, \eqref{slim} implies (see, e.g.,
Theorem VIII.24 in \cite{Reed:Simon})
\begin{equation}\label{limits:2}
\begin{split}
\slim_{n\rightarrow\infty}
\EE_{\mathbf{A}_n+\mathbf{V}_n}((-\infty, 0))
 & = \EE_{\bB}((-\infty, 0)), \\ \slim_{n\rightarrow\infty}
\EE_{\mathbf{A}_n+\mathbf{V}_n}(( 0, \infty))
 &=
\EE_{\bB}(( 0, \infty)),
\end{split}
\end{equation}
and hence
\begin{equation}\label{limits:4}
\slim_{n\rightarrow\infty}
\EE_{\mathbf{A}_n+\mathbf{V}_n}(\{0\})
 =
\EE_{\bB}(\{0\})=0.
\end{equation}

Let $\widehat{\mathbf{A}}_n$ and $\widehat{\mathbf{V}}_n$
denote the parts of the operators $\mathbf{A}_n$ and
$\mathbf{V}_n$ associated with their invariant finite
dimensional subspace $\widehat \fH^{(n)}=\fH_0^{(n)}\oplus
\fH_1^{(n)}$, where $\fH_0^{(n)}=\Ran \cP_n^{(0)}$ and
$\fH_1^{(n)}=\Ran \cP_n^{(1)}$. By Theorem \ref{n0wings} the
subspace (of the finite dimensional Hilbert space $\widehat
\fH_n$)
\begin{equation*}
\Ran \EE_{{\widehat{\mathbf{A}}}_n +
\widehat{\mathbf{V}}_n}((-\infty,0)) \oplus
\big(\Ker(\widehat{\mathbf{A}}_n+ \widehat{\mathbf{V}}_n )\cap
\fH_0^{(n)} \big) \subset \fH^{(n)}
\end{equation*}
is a graph subspace
\begin{equation*}
\cG\big(\Ran \EE_{\widehat{\mathbf{A}}_n}((-\infty,0))\oplus
(\Ker(\widehat{\mathbf{A}}_n )\cap \fH_0^{(n)}), X_n \big)
\end{equation*}
for some $X_n\in \cB(\fH_0^{(n)},\fH_1^{(n)})$, $n\in \N$.
Since $X_n$ is of finite rank, $\|X_n\|\leq 1$ by \eqref{net}.

Applying Theorem \ref{frpr} in the Appendix one arrives at the
inequality
\begin{equation}\label{nunu}
\|\EE_{\widehat{\mathbf{A}}_n +
\widehat{\mathbf{V}}_n}((-\infty,0))\oplus \widehat \cS^{(n)}
- \EE_{\widehat{\mathbf{A}}_n}((-\infty,0))\oplus \widehat
\cT^{(n)} \| = \frac{\|X_n\|}{\sqrt{1+\|X_n\|^2}}\le
\frac{\sqrt{2}}{2},
\end{equation}
where $\widehat \cS^{(n)}$ and $\widehat \cT^{(n)}$ are  the
orthogonal projections in $\fH^{(n)}$ onto the subspaces
$\Ker(\widehat{\mathbf{A}}_n+ \widehat{\mathbf{V}}_n  )\cap
\fH_0^{(n)} $ and  $\Ker(\widehat{\mathbf{A}}_n )\cap
\fH_0^{(n)}$, respectively.

The subspaces $\Ran( \widehat \cS^{(n)})$ and $\Ran( \widehat
\cT^{(n)})$ of the space $\fH^{(n)}$ are naturally imbedded
into the total Hilbert space $\fH$. Denoting by $\cS^{(n)}$
and $\cT^{(n)}$ the corresponding orthogonal projections in
$\fH$ onto these subspaces \eqref{nunu} yields the estimate
\begin{equation}\label{xoroscho}
\|\EE_{\mathbf{A}_n + \mathbf{V}_n}((-\infty,0))\oplus
\cS^{(n)}- \EE_{\mathbf{A}_n}((-\infty,0))\oplus \cT^{(n)}\|
\le \frac{\sqrt{2}}{2} .
\end{equation}

{}From \eqref{limits:0} it follows that
\begin{equation}\label{x1}
\slim_{n\to \infty} \cT^{(n)}=0.
\end{equation}
Meanwhile, by \eqref{limits:4}
\begin{equation}\label{x2}
\slim_{n\to \infty} \cS^{(n)}=0.
\end{equation}

Combining \eqref{xoroscho} -- \eqref{x2} and passing to the limit $n\to
\infty$, by the lower semicontinuity of the spectrum (see, e.g.,
\cite{Kato}, Sec.\ VIII.1.2) one concludes that
\begin{equation*}
\|P-Q\|=\frac{\|X\|}{\sqrt{1+\|X\|^2}}\leq \frac{\sqrt{2}}{2},
\end{equation*}
where $Q$ is the orthogonal projection in $\fH$ onto the subspace $\fQ$
\eqref{Q:Def}. This proves that the operator $X$ is a contraction. The
proof of (i) under the additional assumption that $\Ker \bB=\{0\}$ is
complete.

Step 2.  Assume now that $\Ker \bB$ is not necessarily
trivial. {}From Theorem \ref{n0wings} it follows that the
subspace $\Ker \bB$ is $\bA$-invariant. Denote by $\widehat
\bA$ and $\widehat \bB$ the corresponding parts of the
operators $\bA$ and $\bB$ associated with the reducing
subspace $\widehat{\fH}=\Ran \EE_{\bB}(\R \setminus\{0\})$.
Clearly, the operator $\widehat \bB$ is an off-diagonal
perturbation of the diagonal operator matrix $\widehat \bA$
with respect to the decomposition $\widehat{\fH} =
\widehat{\fH}_0 \oplus \widehat{\fH}_1$, where
\begin{equation*}
\begin{split}
\widehat{\fH}_0 &:= \fH_0\ominus \big(\Ran
\EE_{\bB}(\{0\})\cap
\fH_0\bigr),\\
\widehat{\fH}_1 &:= \fH_1\ominus \big(\Ran
\EE_{\bB}(\{0\})\cap \fH_1\big),
\end{split}
\end{equation*}
and $\Ker\widehat{\bB} = \{0\}$. Moreover, Hypothesis
\ref{diag} is satisfied with the replacements
$\fH\longrightarrow \widehat{\fH}$,
$\fH_0\longrightarrow\widehat{\fH}_0$,
$\fH_1\longrightarrow\widehat{\fH}_1$, and
$\bA\longrightarrow\widehat{\bA}$,
$\bB\longrightarrow\widehat{\bB}$.

By the first part of the proof the subspace
$\Ran\EE_{\bB}\bigl((-\infty,0)\bigr)$, naturally imbedded
into the Hilbert space $\widehat{\fH}$, is the  graph of a
contraction $\widehat X$,
\begin{equation}\label{deso}
\widehat X:\ \widehat{\fH}_0\to \widehat{\fH}_1.
\end{equation}

Clearly, the $\bB$-invariant subspace
\begin{equation*}
\Ran\EE_{\bB}\bigl((-\infty,0)\bigr)\oplus
\big(\Ran\EE_{\bB}(\{0\})\cap \fH_0\big)
\end{equation*}
is a graph subspace $\cG(\fH_0, X)$ associated with the
subspace $\fH_0$ in decomposition \eqref{decom} where the
operator $X$ is given by
\begin{equation}\label{des}
 X f=\begin{cases}
\widehat X f& \text{if}\quad f\in
\Ran\EE_{\bB}((-\infty,0))\,\, (\text{naturally imbedded into
}\widehat{\fH} ),
\\
0& \text{if}\quad f\in \Ran\EE_{\bB}(\{0\})\cap \fH_0= \Ker
A_0\cap\Ker V^\ast.
\end{cases}
\end{equation}
Since $\widehat X$ is a contraction by hypothesis, the
operator $X$ is also a contraction satisfying the properties
\eqref{leer} (with $\lambda=0$):
\begin{equation*}
\Ker A_0\cap\Ker V^\ast\subset \Ker X \text{ \, and \, } \Ran
\EE_{\bB}(\{0\})\cap\fH_1=\Ker A_1\cap\Ker V\subset\Ker X^\ast
\end{equation*}
using \eqref{deso} and \eqref{des}. The proof of (i) is
complete.

(ii) If at least one of  the subspaces $\Ker A_0$ or $\Ker
A_1$ is trivial, then (in the notations above) we have that
\begin{equation*}
(Uf, A_1 Uf)-(f, A_0 f)>0,
\end{equation*}
and, therefore,  equality \eqref{ravno} cannot be satisfied
for $\mu=1$. Hence $\mu=1$ is not a singular number of the
contraction $X$ which proves that the operator $X$ is a strict
contraction, i.e.,
\begin{equation*}
\|Xf\| < \|f\| , \quad\text{for any}\quad f\neq 0.
\end{equation*}
The proof of (ii) is complete.

(iii)  Under   Hypotheses \ref{diag} (with $\lambda=0$) the
fact that the spectra of the operators $A_0$ and $A_1$ are
separated, i.e., $d=\dist\{\spec(A_0), \spec(A_1)\}>0$, means
that at least one of the subspaces $\Ker A_0$  and $\Ker A_1$
is trivial. Therefore, the following estimate holds
\begin{equation}\label{nere}
d\|f\|^2<\big((Uf, A_1 Uf)-(f, A_0 f)\big)
\end{equation}
and, hence,  from \eqref{ravno} one derives the inequality
\begin{equation*}
d<\frac{\mu^2-1}{\mu} \Re (V^\ast f, Uf)\le\frac{1-\mu^2}{\mu}
\|V\|,
\end{equation*}
which proves that the operator $X$ does not have singular values outside
the interval $[0,\nu]$,  where
\begin{equation*}
\nu=\tan \bigg(\frac{1}{2}\arctan \frac{2\|V\|}{d} \bigg)<1.
\end{equation*}
Using the same strategy as in the  proof of (ii) one arrives
to the conclusion that $X$ is a uniform contraction satisfying
the norm estimate
\begin{equation*}
  \|X\|\le \tan \bigg ( \frac{1}{2}\arctan \frac{2\|V\|}{d} \bigg )<1,
\end{equation*}
which proves the upper bound \eqref{normasnizu}.
\end{proof}

\begin{remark}\label{ri}
The operator $X$ referred to in Theorem \ref{Davis:Kahan} is a
contractive solution to the Riccati equation
\begin{equation}\label{Riccati:first}
A_1 X - X A_0 - X V X + V^\ast = 0
\end{equation}
with the property that
\begin{equation}
\spec(A_0 + VX)\subset (-\infty,\lambda],
\end{equation}
since the operator  $A_0+V X$ is similar to the part of $\bB$ associated
with the subspace $\fQ$ (see, e.g., \cite{Albeverio}) and $\sup (\spec
(\bB|_\fQ))\leq \lambda$ by definition \eqref{Q:Def} of the invariant
subspace $\fQ$. The similarity of $A_0+V X$ and $\bB|_\fQ$ can also be seen
directly  from the identity
\begin{equation*}
(x+Xx, \bB(x+Xx)) =(x, (A_0 + X^\ast V^\ast)(I+X^\ast X)x)=(x,
(I+X^\ast X) (A_0 + V X)x)
\end{equation*}
valid for any $x\in\fH_0$.
\end{remark}

The result of Theorem \ref{frpr} in the Appendix shows that Theorem
\ref{Davis:Kahan} admits the following  equivalent formulation
 in terms of the corresponding spectral projections.

\begin{theorem}\label{thm:PQ}
Assume Hypothesis \ref{diag}. Denote by $Q$ the orthogonal
projection in $\fH$ onto the subspace $\fQ$ \eqref{Q:Def}
 and by $P$ the orthogonal projection onto
$\fH_0$. Then:

\emph{(i)} $\|P-Q\|\leq \sqrt{2}/2$,

\emph{(ii)} If  either $\Ker (A_0-\lambda)=\{0\}$ or $\Ker
(A_1-\lambda)=\{0\}$, then
\begin{equation*}
\pm \frac{\sqrt{2}}{2}\notin \specpp(P-Q).
\end{equation*}

\emph{(iii)} If
\begin{equation}\label{razdelenu2}
d=\dist(\spec(A_0), \spec (A_1))>0,
\end{equation}
then
\begin{equation}\label{normasnizu2}
\|P-Q\|\le \sin  \left(\frac{1}{2}\arctan
\frac{2\|V\|}{d}\right) < \frac{\sqrt{2}}{2}.
\end{equation}
\end{theorem}

\begin{remark}\label{gapgap}
Note that if \eqref{razdelenu2} holds, then the interval
$(\sup \spec(A_0), \inf \spec (A_1))$ belongs to the resolvent
set of the operator $\bB$. Although this fact is well known
(see \cite{Adamyan:Langer:95}, \cite{Davis:Kahan}) we present
a particularly simple and short alternative proof.

Given $\lambda\in (\sup \spec(A_0), \inf \spec (A_1))$, the
subspace $\Ran \EE_{\bB}((-\infty, \lambda))$ is a graph subspace
$\cG(\fH_0, X_\lambda) $ where $X_\lambda$ is a strictly
contractive solution to the Riccati equation
\eqref{Riccati:first}. By  a uniqueness result (see Corollary 6.4
(i) in \cite{Kostrykin:Makarov:Motovilov:2}) the solution
$X_\lambda$ does not depend on
\begin{equation*}
\lambda\in (\sup \spec(A_0), \inf \spec (A_1)).
\end{equation*}
Therefore, $\EE_{\bB}\big((\sup\spec(A_0), \inf \spec
(A_1))\big)=0$ which proves the claim.
\end{remark}

As a by-product of our considerations we also get the
following important properties of the subspaces $\Ker
(I_{\fH_0}-X^\ast X)$ and $\Ker (I_{\fH_1}-X X^\ast)$. They
will be used in Sections \ref{sec:unique} and
 \ref{subspaces} below.

\begin{lemma}\label{side}
Let $X$  be the operator referred to in Theorem
\ref{Davis:Kahan}. Then
\begin{equation}\label{odno}
\Ker (I_{\fH_0}-X^\ast X)\subset \Ker (A_0-\lambda),
\end{equation}
\begin{equation}\label{drugoe}
\Ran X|_{\Ker (I_{\fH_0}-X^\ast X)} \subset \Ker
(A_1-\lambda),
\end{equation}
and
\begin{equation}\label{bz}
\Ker (I_{\fH_0}-X^\ast X)\subset \Ker(XVX-V^\ast).
\end{equation}
Moreover, the subspace $\Ker (I_{\fH_0}-X^\ast X)$ reduces
both the operators $V X$ and $V V^\ast$.

Similarly,
\begin{equation}\label{odnoo}
\Ker (I_{\fH_1}-X X^\ast)\subset \Ker (A_1-\lambda),
\end{equation}
\begin{equation}\label{drugoee}
\Ran X^\ast|_{\Ker (I_{\fH_1}-X X^\ast)}\subset \Ker
(A_0-\lambda),
\end{equation}
\begin{equation}\label{bz1}
\Ker (I_{\fH_1}-X X^\ast)\subset \Ker(X^\ast V^\ast X^\ast -
V),
\end{equation}
and the subspace $\Ker (I_{\fH_1}-X X^\ast)$ reduces both the operators
$X^\ast V^\ast$ and $V^\ast V$.
\end{lemma}

\begin{proof}
Let $0\ne f\in \Ker (I_{\fH_0}-X^\ast X)$, that is, $|X|f=f$, where
$|X|=(X^\ast X)^{1/2}$. Then the elements $f\oplus Uf$ and $\bB((-f)\oplus
Uf)$ are orthogonal, where $U$ is the partial isometry from the polar
decomposition $X=U|X|$ (see the proof of Theorem \ref{Davis:Kahan} part
(i)). This means that
\begin{equation}
(Uf, (A_1-\lambda) Uf)=(f, (A_0-\lambda)f),
\end{equation}
which implies \eqref{odno} and \eqref{drugoe} (under
Hypothesis \ref{diag}).

In order to prove \eqref{bz} and that the subspace
$\Ker(I_{\fH_0} - X^\ast X)$ is both $VX$- and $V
V^\ast$-invariant we proceed as follows.

By Remark \ref{ri} the operator $X$ solves the Riccati
equation
\begin{equation}\label{ri0}
A_1 X - X A_0 - X V X + V^\ast = 0
\end{equation}
and, hence,
\begin{equation}\label{ri1}
X^\ast A_1  -  A_0 X^\ast - X^\ast V^\ast X^\ast + V = 0,
\end{equation}
which in particular implies that
\begin{equation}\label{ri2}
X^\ast A_1 X  -  A_0 X^\ast X - X^\ast V^\ast X^\ast X + V X =
0.
\end{equation}
For any  $f\in \Ker (I_{\fH_0}-X^\ast X)$ inclusions
\eqref{odno} and \eqref{drugoe} yield
\begin{equation*}
(A_1 X  -  X A_0)f = (X^\ast A_1 X  - A_0) f = (X^\ast A_1 X
-  A_0 X^\ast X)f = 0.
\end{equation*}
Thus, from \eqref{ri0} and \eqref{ri2} it follows that
\begin{equation}\label{comb1}
V^\ast f = X V X f, \qquad f\in \Ker (I_{\fH_0}-X^\ast X),
\end{equation}
which proves \eqref{bz} and the representation
\begin{equation}\label{comb2}
VXf = X^\ast V^\ast X^\ast X f=X^\ast V^\ast f, \qquad f\in
\Ker (I_{\fH_0}-X^\ast X).
\end{equation}

Combining \eqref{comb1} and \eqref{comb2} proves that
\begin{equation*}
(I_{\fH_0} -X^\ast X) V X f = V X f - X^\ast X V X f = V X f -
X^\ast V^\ast f = 0
\end{equation*}
for any $f\in \Ker (I_{\fH_0}-X^\ast X)$. That is, the
subspace $\Ker(I_{\fH_0} - X^\ast X)$ is $VX$-invariant. From
\eqref{comb2} it follows  that $\Ker(I_{\fH_0} - X^\ast X)$ is
also $X^\ast V^\ast $-invariant and, hence, $\Ker(I_{\fH_0} -
X^\ast X)$ reduces the operator $VX$.

Equality \eqref{comb1} implies that
\begin{equation}
VV^\ast f = V X V X f, \qquad f\in \Ker (I_{\fH_0}-X^\ast X)
\end{equation}
which proves that $\Ker (I_{\fH_0}-X^\ast X)$ is $V
V^\ast$-invariant, since $\Ker (I_{\fH_0}-X^\ast X)$ is
already proven to be a $V X$-invariant subspace. Since $V
V^\ast$ is self-adjoint, the subspace $\Ker (I_{\fH_0}-X^\ast
X)$ reduces $V V^\ast$. The proof of \eqref{odnoo},
\eqref{drugoee}, and \eqref{bz1} is similar.
\end{proof}

\begin{remark}
It follows from Lemma \ref{side} that the multiplicity $m$ of
the singular value $\mu=1$ of the operator $X$ satisfies the
inequality
\begin{equation}\label{kratnost}
m \leq  \min \{ \dim \Ker(A_0-\lambda),\dim
\Ker(A_1-\lambda)\}.
\end{equation}
Equivalently,
\begin{equation}\label{ineq:kerne:1}
\begin{split}
\dim\Ker \left(Q-P-\frac{\sqrt{2}}{2}\right) & =
\dim\Ker\left(Q-P+\frac{\sqrt{2}}{2}\right)\\ &\leq  \min
\{\dim \Ker(A_0-\lambda),\dim \Ker(A_1-\lambda)\},
\end{split}
\end{equation}
where $P$ and $Q$ are orthogonal projections in $\fH$ onto the subspaces
$\fH_0$ and $\fQ$ \eqref{Q:Def}, respectively. The subspaces $\Ker
(I_{\fH_0}-XX^\ast)$ and $\Ker (I_{\fH_1}-X^\ast X)$ will be studied in
Section \ref{subspaces} below.
\end{remark}

\section{Lower Bound}\label{sec:lower:bound}

In this section we derive the lower bound on the norm of the difference of
the orthogonal projections onto the $\bA$-invariant subspace $\fH_0$ and
the $\bB$-invariant subspace $\fQ$ given by \eqref{Q:Def}.

\begin{theorem}\label{XXL}
Assume Hypothesis \ref{diag}. Let ${\delta_-}$ (${\delta_+}$) denote the shift of the bottom (top,
respectively) of the spectrum of the operator $\mathbf{A}$ under the perturbation $\mathbf{V}$,
i.e.,
\begin{equation*}
{\delta_-} = \inf\spec(\mathbf{A}) - \inf\spec(\bB),\qquad {\delta_+} = \sup\spec(\bB) -
\sup\spec(\mathbf{A}).
\end{equation*}
Then the solution $X$  to the Riccati equation referred to in Theorem \ref{Davis:Kahan} satisfies
the lower bound
\begin{equation}\label{X:lower}
\|X\|\geq  \frac{\delta}{\|V\|},
\end{equation}
where
\begin{equation*}
\delta = \max\{\delta_-, \delta_+\}\leq \|V\|.
\end{equation*}
Equivalently,
\begin{equation}\label{lower}
\|Q-P\| \geq \frac{{\delta}}{\sqrt{\delta^2 + \|V\|^2}},
\end{equation}
where $P$ and $Q$ are orthogonal projections onto the subspace $\fH_0$ and
$\fQ$ \eqref{Q:Def}, respectively.
\end{theorem}

\begin{remark}
{}From a general perturbation theory for off-diagonal perturbations it follows that $\delta_-\geq
0$ and $\delta_+\geq 0$. For the proof of this fact we refer to
\cite{Kostrykin:Makarov:Motovilov:3}.
\end{remark}

\begin{proof}[Proof of Theorem \ref{XXL}]
{}From Theorem \ref{Davis:Kahan}  it follows that
$\Ran\EE_{\bB}((-\infty,\lambda))\oplus(\fN_0\cap\Ker V^\ast)$ is
the  graph subspace $\cG(\fH_0, X)$
associated with the subspace $\fH_0$ in the decomposition
$\fH=\fH_0 \oplus \fH_1$ and by Remark \ref{ri}
 the operator $X$ solves the Riccati equation
\eqref{Riccati:first}. It is well known (see, e.g., \cite{Albeverio}) that
in this case
\begin{equation}
\spec (\bB)=\spec(A_0+VX)\cup \spec(A_1-X^\ast V^\ast)
\end{equation}
and, hence,
\begin{equation}\label{5.37}
\inf\spec(\bB)\geq \inf\spec(\mathbf{A}) -\|V\| \|X\|
\end{equation}
and
\begin{equation}\label{5.38}
\sup\spec(\bB)\leq \sup\spec(\mathbf{A}) +\|V\| \|X\|,
\end{equation}
which  proves the lower bounds \eqref{X:lower} and \eqref{lower}
using Theorem \ref{frpr} in the Appendix.
\end{proof}

\begin{remark}\label{lower:sharp}
Let $\fH_0=\fH_1=\C$, $A_0$ and $A_1$ are reals with $A_1 >
A_0$, and $V\in \C$. Then
\begin{equation*}
\begin{split}
& \spec \begin{pmatrix} A_0 & V \\ V^\ast & A_1 \end{pmatrix}  \\
& \quad = \bigg \{A_0 -\|V\| \tan\left(\frac{1}{2}\arctan
\frac{2\|V\|}{ A_1- A_0}\right), \, A_1 +\|V\|
\tan\left(\frac{1}{2}\arctan \frac{2\|V\|}{ A_1-
A_0}\right)\bigg \},
\end{split}
\end{equation*}
which  can easily be seen by solving the characteristic equation
\begin{equation*}
(A_0-\lambda)(A_1-\lambda)-\|V\|^2=0.
\end{equation*}
The upper bounds \eqref{normasnizu} and \eqref{normasnizu2} give
the exact value of the norms $\|X\|$ and $\|P-Q\|$, respectively,
and, hence,  the bounds \eqref{normasnizu} and \eqref{normasnizu2}
are sharp. In this case
\begin{equation*}
{\delta_-} = {\delta_+} = \|V\| \tan\left(\frac{1}{2}\arctan \frac{2 \|V\|}{A_1 - A_0} \right)
\end{equation*}
and
\begin{equation*}
\|P-Q\| = \sin\left(\frac{1}{2}\arctan \frac{2 \|V\|}{A_1 -
A_0}, \right)
\end{equation*}
which shows that estimates \eqref{lower} and \eqref{X:lower} are also sharp.
\end{remark}

\section{Riccati Equation: Uniqueness}\label{sec:unique}

Under Hypothesis \ref{diag} Theorem \ref{Davis:Kahan} and
Remark \ref{ri} guarantee the existence of a contractive
solution to the Riccati equation
\begin{equation}
\label{wieder:Riccati} A_1 X - X A_0 - XVX + V^\ast = 0
\end{equation}
with the properties
\begin{equation}\label{i-ii}
\begin{split}
&\textrm{(i)}\quad \Ker(A_0-\lambda)
 \cap \Ker V^\ast \subset \Ker X,\\
  &\textrm{(ii)}\quad
\spec(A_0 + VX)\subset (-\infty,\lambda].
\end{split}
\end{equation}
If, in addition,
\begin{equation*}
d=\dist (\spec (A_0) , \spec(A_1)) > 0,
\end{equation*}
then $X$ is a unique contractive solution to the
Riccati equation (see Corollary 6.4 (i) in
\cite{Kostrykin:Makarov:Motovilov:2}; cf.\
\cite{Adamyan:Langer:Tretter:2000a}).

The following uniqueness result shows that there is no other
solution to the Riccati equation \eqref{wieder:Riccati} with
properties \eqref{i-ii}.

\begin{theorem}\label{Clubshilfe}
Assume Hypothesis \ref{diag}. Then a contractive solution to
the Riccati equation \eqref{wieder:Riccati} satisfying the
properties \eqref{i-ii}  is unique.
\end{theorem}

\begin{proof}
Exactly the same reasoning as in the  proof of Theorem
\ref{Davis:Kahan} (i) allows to conclude that without loss of
generality one may assume that $\Ker (\bB-\lambda)=0$.

In this case the spectral subspace $\Ran \EE_{\bB}((-\infty,
\lambda))$  is the graph of a contraction $X$ which solves
\eqref{wieder:Riccati}. Assume
that $Y\in\cB(\fH_0,\fH_1)$ is a (not necessarily contractive)
solution to \eqref{wieder:Riccati} different from $X$. Since the
graph of $X$ is a spectral subspace of $\bB$, one concludes that
the orthogonal projections onto the graphs $\cG(\fH_0,X)$ and
$\cG(\fH_0,Y)$ of the operators $X$ and $Y$ commute.

We claim that $\cG(\fH_0,X)^\perp\cap\cG(\fH_0,Y)$ is
nontrivial. To show this we set
\begin{equation*}
\fL:= \fH_0\ominus \Ran P|_{\cG(\fH_0,X)\cap\cG(\fH_0,Y)},
\end{equation*}
where $P$ is the orthogonal projection in $\fH$ onto the subspace
$\fH_0$. Any $z\in\cG(\fH_0,X)\cap\cG(\fH_0,Y)$ admits the
representations
\begin{equation*}
\begin{split}
z= x\oplus Xx\quad \text{for some}\quad x\in\fL^\perp,\\
z= y\oplus Yy\quad \text{for some}\quad y\in\fL^\perp,
\end{split}
\end{equation*}
where $\fL^\perp:=\fH_0\ominus\fL$. Obviously, $x = y = Pz$,
and, therefore, $XPz=YPz$ for any
$z\in\cG(\fH_0,X)\cap\cG(\fH_0,Y)$. Thus,
$Y|_{\fL^\perp}=X|_{\fL^\perp}$  and
\begin{equation*}
\cG(\fH_0,X)\cap\cG(\fH_0,Y) = \{x+Xx|\ x\in\fL^\perp\}.
\end{equation*}
Hence,  $X\neq Y$ if and only if the subspace $\fL$ is
nontrivial.

Note that
\begin{equation*}
(x_0\oplus Y x_0, x\oplus Y x) =0,\qquad x:= (I_{\fH_0}+
Y^\ast Y)^{-1} y
\end{equation*}
for any $x_0\in\fL^\perp$ and any $y\in\fL$. Hence,
\begin{equation}\label{bei}
\cG(\fH_0,Y) \ominus (\cG(\fH_0,X)\cap\cG(\fH_0,Y))= \{x+Yx|\
x=(I_{\fH_0}+Y^\ast Y)^{-1} y,\, y\in\fL\}.
\end{equation}
Since the orthogonal projections onto the graph subspaces
$\cG(\fH_0,X)$ and $\cG(\fH_0,Y)$ commute, by \eqref{bei} we conclude
that the subspace
\begin{equation*}
\cG(\fH_0,X)^\perp\cap\cG(\fH_0,Y) = \cG(\fH_0,Y)\ominus
(\cG(\fH_0,X)\cap\cG(\fH_0,Y))
\end{equation*}
is nontrivial.

For any $z\in\fH$, $z\neq 0$ such that $z\in\cG(\fH_0, Y)$ and
$z\perp\cG(\fH_0, X)$ we have
\begin{equation}\label{varvar}
(z, \bB z)>\lambda.
\end{equation}
Therefore, for the operator $A_0+VY$ the condition (ii) does
not hold, since the spectrum of $A_0+VY$ coincides with that
of the restriction of $\bB$ onto its invariant subspace
$\cG(\fH_0, Y)$ and by \eqref{varvar} the operator
$\bB|_{\cG(\fH_0, Y)}$ has points of  the spectrum to the
right of the point $\lambda$.  The proof is complete.
\end{proof}

\begin{remark}\label{remarka:2}
If $\lambda\in \R$ is a multiple eigenvalue of the operator
$\bB$ and both $\Ker(\bB-\lambda)\cap \fH_0\ne \{0\}$ and
$\Ker(\bB-\lambda)\cap \fH_1\ne \{0\}$, it follows from
Theorem \ref{n0wings} that the Riccati equation
\eqref{wieder:Riccati} has uncountably many bounded solutions
(even if $\fH$ is finite-dimensional). This can also be seen
directly. Let
\begin{equation*}
T: \Ker(A_0-\lambda)\cap\Ker V^\ast
\longrightarrow\Ker(A_1-\lambda)\cap\Ker V
\end{equation*}
be any bounded operator acting from $\Ker(A_0-\lambda)\cap\Ker
V^\ast$ to $\Ker(A_1-\lambda)\cap\Ker V$. The bounded operator
$X\in\cB(\fH_0,\fH_1)$ defined by
\begin{equation}\label{X:Def}
\widetilde X f=\begin{cases}
T f& \text{ if }\, f\in \Ker(A_0-\lambda)\cap\Ker V^\ast\\
0& \text{ if }\, f\in \Ker(A_0-\lambda)\ominus
(\Ker(A_0-\lambda)\cap\Ker V^\ast)
\end{cases}
\end{equation}
satisfies the equation
\begin{equation}
\big((A_1-\lambda I_{\fH_1})\widetilde X - \widetilde X (
A_0-\lambda I_{\fH_0}) - \widetilde X V\widetilde  X +
V^\ast\big)f=0 \quad \text{for all}\quad f\in \fH_0,
\end{equation}
and, thus, it is also a solution to \eqref{wieder:Riccati}. If
$\dim\Ker(A_0-\lambda)\cap\Ker V^\ast=\infty$,
$\dim\Ker(A_1-\lambda)\cap\Ker V=\infty$, and $T$ is a closed
densely defined unbounded operator from
$\Ker(A_0-\lambda)\cap \Ker V^\ast$ to  $\Ker(A_1-\lambda)\cap
\Ker V$, then the operator $\widetilde X$ defined by
\eqref{X:Def} is an unbounded solution to the Riccati equation
\eqref{wieder:Riccati} in the sense of Definition 4.2 in
\cite{Kostrykin:Makarov:Motovilov:2}.
\end{remark}

Our next result is the following uniqueness criterion.

\begin{theorem}\label{4.11a}
Assume Hypothesis \ref{diag}. A contractive solution $X$ to
the Riccati equation \eqref{wieder:Riccati} is a unique
contractive solution if and only if the graph of $X$ is a
spectral subspace of the operator $\bB$ and $\mu=1$ is not an
eigenvalue of the operator $|X|$, the absolute value of $X$.
In this case, $X$ is a strict contraction and
\begin{equation*}
\text{either}\quad \EE_\bB((-\infty, \lambda))=\cG(\fH_0,
X)\quad\text{or}\quad \EE_\bB((\lambda,+\infty))=\cG(\fH_1,
-X^\ast).
\end{equation*}
\end{theorem}
\begin{proof}
``If Part". Since $\mu=1$ is not an eigenvalue of the
contraction $X$, it follows that $X$ is a strict contraction.
Then the claim follows from Corollary 6.4 of
\cite{Kostrykin:Makarov:Motovilov:2}

``Only If Part". Assume that $X$ is the  unique contractive
solution to the Riccati equation.  Then $X$ coincides with the
operator referred to in Theorem \ref{Davis:Kahan}. We need to
prove that the graph of $X$ is a spectral subspace of $\bB$
and that $\Ker(I_{\fH_0}-X^\ast X)=\{0\}.$

We will prove these statements by reduction to contradiction. If
the graph of $X$ is not a spectral subspace of $\bB$, then by
Remark \ref{remark:on:Q} and Theorem \ref{Davis:Kahan} both
$\Ker(A_0-\lambda)\cap\Ker V^\ast$ and $\Ker(A_1-\lambda)\cap\Ker
V$ are nontrivial. Let $T$ be any contractive operator from
$\Ker(A_0-\lambda)\cap\Ker V^\ast$ to $\Ker(A_1-\lambda)\cap\Ker
V$ with $\Ker T\neq \{0\}$. Then the operator $\widetilde X$
defined by \eqref{X:Def} is also a contractive solution to the
Riccati equation (see Remark \ref{remarka:2}). Since by Theorem
\ref{Davis:Kahan} $\Ker (A_0-\lambda)\cap \Ker V^\ast \subset \Ker
X$, where $X$ is the contractive solution to the Riccati equation
referred to in Theorem \ref{Davis:Kahan}, the contractive solution
$\widetilde X$ to the Riccati equation  is different from $X$ by
construction. A contradiction.

Assume now that $\mu=1$ is an eigenvalue of $|X|$, that is,
$\Ker(I_{\fH_0} - X^\ast X)$ is nontrivial.  By Lemma
\ref{side}
\begin{equation}\label{XVX}
\Ker(I_{\fH_0} - X^\ast X)\subseteq \Ker(XVX-V^\ast)
\end{equation}
and $\Ker(I_{\fH_0} - X^\ast X)$ reduces both $A_0$  and $VX$.
Applying Theorem 6.2 in \cite{Kostrykin:Makarov:Motovilov:2}
we conclude that the Riccati equation \eqref{wieder:Riccati}
has a contractive solution $Y$ such that $\Ker(I_{\fH_0} -
X^\ast X)=\Ker(I_{\fH_0}+Y^\ast X)$. This solution necessarily
differs from $X$ which contradicts the hypothesis that $X$ is
a unique contractive solution.

{}From Theorem \ref{n0wings} it follows now  that the graph of $X$ is the
spectral subspace of the operator~$\bB$ and
\begin{equation*}
\cG(\fH_0, X)=
\begin{cases}
\EE_\bB((-\infty, \lambda)) & \text{if}\quad \Ker(A_0-\lambda)\cap\Ker V^\ast=\{0\}\\
\EE_\bB((-\infty, \lambda]) & \text{if}\quad
\Ker(A_1-\lambda)\cap\Ker V=\{0\}
\end{cases},
\end{equation*}
which proves the remaining statement of the theorem.
\end{proof}

As an immediate consequence of Theorem \ref{4.11a} we get the
following results.

\begin{corollary}\label{cor:strictly:contractive}
Let $X:\fH_0\rightarrow\fH_1$ be the solution to the Riccati
equation referred to in Theorem \ref{Davis:Kahan}. It is the
unique contractive solution if and only if it is strictly
contractive.
\end{corollary}

\begin{remark}
Let $X$ be an arbitrary strictly contractive solution to the
Riccati equation \eqref{wieder:Riccati}.  In general this
solution has not to be the unique contractive solution or to
be isolated point of the set of all solutions (cf.\ Remark
\ref{remarka:2}).
\end{remark}

\begin{corollary}\label{isolated}
Let $X:\fH_0\rightarrow\fH_1$ be the solution to the Riccati
equation referred to in Theorem \ref{Davis:Kahan}. It is an
isolated point (in the operator norm topology) in the set of all
solutions  to
the Riccati equation \eqref{wieder:Riccati} if and only if either
$\Ker(A_0-\lambda)\cap\Ker V^\ast =\{0\}$ or
$\Ker(A_1-\lambda)\cap\Ker V =\{0\}$.
\end{corollary}

\begin{proof}
By Remark \ref{remark:on:Q} and Theorem \ref{Davis:Kahan} the
graph $\cG(\fH_0,X)$ is associated with a spectral subspace of
the operator $\bB$ if and only if either
$\Ker(A_0-\lambda)\cap\Ker V^\ast =\{0\}$ or
$\Ker(A_1-\lambda)\cap\Ker V =\{0\}$. Now the claim
 follows from Theorem 5.3 in \cite{Kostrykin:Makarov:Motovilov:2}.
\end{proof}

\begin{remark}
If either $\Ker(A_0-\lambda)\cap\Ker V^\ast =\{0\}$ or
$\Ker(A_1-\lambda)\cap\Ker V =\{0\}$ holds, Theorem 6.2 in
\cite{Kostrykin:Makarov:Motovilov:2} allows  to construct \emph{all} contractive solutions to
the Riccati equation from that referred to in Theorem
\ref{Davis:Kahan}.
\end{remark}

\section{More on the Subspaces $\Ker (I_{\fH_0}-X^\ast X)$ and
 $\Ker (I_{\fH_1}-XX^\ast)$}\label{subspaces}

The main goal of this section is to prove the  fact that the
subspace $\Ker (I_{\fH_0}-X^\ast X)$ associated with the
operator $X$ referred to in Theorem \ref{Davis:Kahan} admits
an  intrinsic description as the maximal $VV^\ast$-invariant
subspace $\fK_0\subset\fH_0$ with the properties
\begin{equation}\label{mp1}
\begin{split}
&\fK_0\subset
\Ker(A_0-\lambda)\cap \overline{\Ran V},\\
&\Ran V^\ast|_{\fK_0}\subset\Ker(A_1-\lambda)
\end{split}
\end{equation}
(see Theorem \ref{karkar} below). Similarly,  the subspace
$\Ker (I_{\fH_1}-XX^\ast)$ can be characterized as the maximal
$V^\ast V$-invariant subspace with the properties
\begin{equation}\label{mp11}
\begin{split}
&\fK_1\subset \Ker(A_1-\lambda)\cap
\overline{\Ran V^\ast},\\
&\Ran V|_{\fK_1}\subset\Ker(A_0-\lambda).
\end{split}
\end{equation}

We start with the observation that the  maximal subspaces
$\fK_0$ and $\fK_1$  with indicated properties do exist and
admit a constructive description.

\begin{lemma}
The subspaces
\begin{equation}\label{K0}
\begin{split}
\fK_0 := & \mathrm{closure}\big\{ x\in\Ker(A_0-\lambda)\cap
\Ran V \big| \\ & \quad V^\ast (V V^\ast)^n x \in
\Ker(A_1-\lambda), \,\, (V V^\ast)^n x \in \Ker(A_0-\lambda)\
\text{for any}\ n\in\N_0 \big\}
\end{split}
\end{equation}
and
\begin{equation}\label{K1}
\begin{split}
\fK_1 := & \mathrm{closure}\big\{ x\in\Ker(A_1-\lambda)\cap
\Ran V^\ast \big| \\ & \quad V (V^\ast V)^n x \in
\Ker(A_0-\lambda), \,\, (V^\ast V)^n x \in \Ker(A_1-\lambda)\
\text{for any}\ n\in\N_0 \big\}
\end{split}
\end{equation}
reduce the operators $VV^\ast$ and $V^\ast V$, respectively.
Moreover, the subspaces $\fK_0$ and $\fK_1$ satisfy the
properties \eqref{mp1} and \eqref{mp11}, respectively.

The subspace $\fK_0$ is  maximal in the sense that if $\cL_0$ is
any other $VV^\ast$-invariant subspace with the properties
\eqref{mp1}, then $\fL_0\subset \fK_0$. Analogously, the subspace
 $\fK_1$ is  maximal in the sense that if $\cL_1$ is  any
other $V^\ast V$-invariant subspace with the properties
\eqref{mp11}, then $\fL_1\subset \fK_1$.
\end{lemma}

\begin{proof}
Clearly, the subspace $\fK_0$ is invariant under the operator
$VV^\ast$ and, therefore, $\fK_0$ reduces $VV^\ast$, since
$VV^\ast$ is self-adjoint. It follows from \eqref{K0} that
\eqref{mp1} holds.

Now, let $\fL_0$ be  an arbitrary closed subspace of
$\Ker(A_0-\lambda)\cap \overline{\Ran V}$ invariant under
$VV^\ast$ such that $V^\ast\fL\subset\Ker(A_1-\lambda)$. Then
$\fL\subset\fK_0$. Indeed, since $\fL$ is invariant under
$VV^\ast$, we have $(VV^\ast)^n
\fL\subset\fL\subset\Ker(A_0-\lambda)$ for any $n\in\N$.
Hence, $V^\ast (VV^\ast)^n \fL\subset
V^\ast\fL\subset\Ker(A_1-\lambda)$ for any $n\in\N_0$ and one
concludes that $\fL\subset\fK_0$.

The maximality of the subspace $\fK_1$ is proven in a similar
way.
\end{proof}

\begin{lemma}\label{prep:lemma:1}
The subspaces $\fK_0$ and $\fK_1$ satisfy the properties that
\begin{align}
\overline{\Ran V^\ast|_{\fK_0}} = \fK_1, &\qquad
\overline{\Ran V|_{\fK_1}} = \fK_0,\label{Zeile:1}\\
\Ran V^\ast|_{\fH_0\ominus\fK_0} \subset \fH_1\ominus\fK_1,
&\qquad  \Ran V|_{\fH_1\ominus\fK_1} \subset
\fH_0\ominus\fK_0.\label{Zeile:2}
\end{align}
\end{lemma}

\begin{proof}
Equations \eqref{Zeile:1} follow from the explicit description
\eqref{K0} and \eqref{K1} of the subspaces $\fK_0$ and
$\fK_1$, respectively.

Let $x\in\fH_0\ominus\fK_0$ be arbitrary. Choose an arbitrary
$y\in\fK_1$ and consider
\begin{equation*}
(y, V^\ast x) = (V y, x).
\end{equation*}
Since, by \eqref{Zeile:1}, $V y\in\fK_0$ we have $(V y, x)=0$. Thus,
$V^\ast x\in\fH_1\ominus\fK_1$ which proves the first inclusion in
\eqref{Zeile:2}. The second inclusion in \eqref{Zeile:2} is proven
similarly.
\end{proof}

\begin{theorem}\label{karkar}
Assume Hypothesis \ref{diag}. Let $X:\fH_0\rightarrow\fH_1$ be
the solution to the Riccati equation \eqref{wieder:Riccati}
referred to in Theorem \ref{Clubshilfe}. Then
\begin{equation}\label{ininin}
\Ker(I_{\fH_0} - X^\ast X) =
\fK_0\quad\text{and}\quad\Ker(I_{\fH_1} - X X^\ast) = \fK_1.
\end{equation}
Moreover, $\overline{\Ran X|_{\fK_0}} = \fK_1$ and
\begin{equation*}
X|_{\fK_0} = - \widehat{S}, \qquad
\widehat{S}:\fK_0\rightarrow \fK_1,
\end{equation*}
where $\widehat{S}=S|_{\fK_0}$ with $S:\fH_0\rightarrow\fH_1$
being the partial isometry with initial space $\overline{\Ran
V}$ and final space $\overline{\Ran V^\ast}$ defined by the
polar decomposition $V^\ast=S(VV^\ast)^{1/2}$. In particular,
$\fK_1=\cG(\fK_0,\widehat{S})$.
\end{theorem}

\begin{proof}
Without loss of generality we assume that $\lambda=0$.

First, we will prove the inclusion
\begin{equation}
\Ker(I_{\fH_0} - X^\ast X) \subset \fK_0.
\end{equation}
It is sufficient to establish that

(a) $\Ker(I_{\fH_0}-X^\ast X)$ reduces $V V^\ast$,

(b) $\Ker(I_{\fH_0} - X^\ast X)  \subset \overline{\Ran
V^\ast} \cap \Ker A_0$,

(c) $V^\ast \Ker(I_{\fH_0} - X^\ast X) \subset \Ker A_1$.

The statement (a) follows from Lemma \ref{side}.

In order to see that (b) holds note that if $z\in\Ker V^\ast \cap
\Ker A_0$, then $Xz=0$, since $\Ker V^\ast \cap \Ker A_0\subset
\Ker X$ by Theorem \ref{Davis:Kahan}.  Therefore, $\Ker V^\ast
\cap \Ker A_0 \perp \Ker(I_{\fH_0}-X^\ast X)$ since
\begin{equation*}
(z,x) = (z, X^\ast X x) = (X z, X x) = 0 \quad \text{for any}\quad
x\in\Ker(I_{\fH_0}-X^\ast X),
\end{equation*}
which proves (b) taking into account that $\Ker(I_{\fH_0} - X^\ast
X) \subset \Ker A_0$ (see Lemma \ref{side}).

To prove (c) we proceed as follows. If $x\in \Ker(I_{\fH_0}-X^\ast
X)$, then $A_0x=0$, since $\Ker(I_{\fH_0} - X^\ast X) \subset \Ker
A_0$. The same reasoning as in the proof of Lemma \ref{side} shows
that (see \eqref{comb1})
\begin{equation*}
V^\ast x = X V X x\qquad\text{and}\qquad V X = X^\ast V^\ast x.
\end{equation*}
Therefore, $(I_{\fH_1}-XX^\ast)V^\ast x=V^\ast x - X X^\ast
V^\ast x=V^\ast x-X V X x=0$, that is, $V^\ast x\in
\Ker(I_{\fH_1}-X X^\ast)$. By Lemma \ref{side}
$\Ker(I_{\fH_1}-XX^\ast)\subset \Ker A_1$, which completes the
proof of (c).

The inclusion $\Ker(I_{\fH_1}-X X^\ast) \subset \fK_1$ is proven similarly.
Hence, we have established that
\begin{equation}\label{kula}
\Ker(I_{\fH_0} - X^\ast X) \subset
\fK_0\quad\text{and}\quad\Ker(I_{\fH_1} - X X^\ast) \subset
\fK_1.
\end{equation}

Now we turn to the proof of the opposite inclusions
\begin{equation}\label{vklu}
\fK_0\subset \Ker(I_{\fH_0}-X^\ast X) \qquad \text{and}\qquad
\fK_1\subset \Ker(I_{\fH_1}-X X^\ast).
\end{equation}
Clearly, the subspaces $\widehat{\fH}_0=\fH_0 \ominus\fK_0$ and
$\widehat{\fH}_1=\fH_1 \ominus\fK_1$ reduce the operators $A_0$ and $A_1$,
respectively, since $\fK_0\subset\Ker A_0$, $\fK_1\subset\Ker A_1$, and
$A_0$, $A_1$ are self-adjoint.  Denote by
$\widehat{A}_0$ and $\widehat{A}_1$ the corresponding parts of the
operators $A_0$ and $A_1$ associated with these subspaces:
\begin{equation*}
\widehat{A}_0 = A_0|_{\widehat{\fH}_0}\qquad\text{and}\qquad
\widehat{A}_1 = A_1|_{\widehat{\fH}_1}.
\end{equation*}
Since by Lemma \ref{prep:lemma:1}  $\Ran V|_{\widehat{\fH}_1}\subset
\widehat{\fH}_0$, the restriction $\widehat{V}$ of the operator $V$ onto
$\fH_1$ is a map from $\widehat{\fH}_1$ to $\widehat{\fH}_0$. By Theorem
\ref{Clubshilfe} the Riccati equation
\begin{equation}\label{urez}
\widehat{A}_1 \widehat{X} - \widehat{X} \widehat{A}_0 -
\widehat{X}\widehat{V}\widehat{X} + \widehat{V}^\ast=0
\end{equation}
has a unique solution $\widehat{X}$ satisfying $\Ker
\widehat{A}_0 \cap \Ker \widehat{V}^\ast\subset
\Ker\widehat{X}$ and
$\spec(\widehat{A}_0+\widehat{V}\widehat{X})\subset(-\infty,0]$.

Let $S:\fH_0\rightarrow\fH_1$ be a partial isometry with the
initial subspace $\overline{\Ran V}$ and final subspace
$\overline{\Ran V^\ast}$ defined by the polar decomposition
$V^\ast = S (VV^\ast)^{1/2}$. {}From Lemma \ref{prep:lemma:1}
it follows that $\Ran S|_{\fK_0}=\fK_1$.

Let
$\widehat{S}:\fK_0\rightarrow\fK_1$ be the restriction of $S$ onto
$\fK_0$: $\widehat{S}=S|_{\fK_0}$. Define  the operator $Y:\
\fH_0\rightarrow\fH_1$ by the following rule
\begin{equation*}
Yx=\begin{cases} \widehat{X}x, & x\in\widehat{\fH}_0\\
-\widehat{S}x, & x\in \fK_0
\end{cases}.
\end{equation*}
Since $\widehat{S}$ maps $\fK_0$ onto $\fK_1$ isometrically,
 one immediately concludes that
\begin{equation}\label{tak}
\fK_0\subset \Ker(I-Y^\ast
Y)\qquad\text{and}\qquad\fK_1\subset \Ker(I-Y Y^\ast).
\end{equation}

We  claim  that the operator $Y$
solves the Riccati equation
\begin{equation}\label{kkk}
A_1 Y - Y A_0 - YV Y + V^\ast =0.
\end{equation}
Indeed, if $x\in\widehat{\fH}_0$ then
\begin{equation}\label{riri}
                (A_1 Y - Y A_0 - Y V Y + V^\ast) x = 0
\end{equation}
as a consequence of \eqref{urez}. If  $x\in\fK_0$,  then  $A_0
x = 0$ (recall that we assumed that  $\lambda = 0$). Moreover,
$\widehat{S} x \in \fK_1$ and hence  $A_1 \widehat{S} x = 0$
resulting in
\begin{equation*}
                (A_1 Y - Y A_0 - Y V Y + V^\ast) x
      = (- Y V Y + V^\ast) x = (-\widehat{S} V \widehat{S} + V^\ast)x=0,
\end{equation*}
where we have used the fact that
$V\widehat{S}x\in\fK_0$ and the equality
\begin{equation*}
\widehat{S} V \widehat{S} x = \widehat{S} (VV^\ast)^{1/2}
\widehat{S}^\ast \widehat{S} x
                = \widehat{S} (VV^\ast)^{1/2} x = V^\ast x.
        \end{equation*}
Therefore, $Y$ solves the Riccati equation \eqref{riri}.

Our next claim is that
\begin{equation}\label{mmm}
\Ker A_0 \cap \Ker V^\ast \subset \Ker Y.
\end{equation}
Since $\fK_0\subset \Ker A_0$, by Lemma \ref{prep:lemma:1} one
concludes that the subspace $\Ker\widehat{A}_0 \cap \Ker\widehat
V^\ast$, naturally imbedded into $\fH_0$,  coincides with  $\Ker A_0 \cap \Ker V^\ast$. One also
concludes that the subspace $\Ker \widehat X$ naturally imbedded
into $\fH_0$ coincides with $\Ker Y$ by the definition of the
operator $Y$. Therefore, \eqref{mmm} follows from the inclusion $
\Ker \widehat A_0 \cap \Ker \widehat V^\ast \subset \Ker \widehat
X$, proving \eqref{mmm}.

Finally, observe that
\begin{equation*}
VS|_{\fK_0} = (VV^\ast)^{1/2} S^\ast S|_{\fK_0} =
(VV^\ast)^{1/2}|_{\fK_0} \geq 0
\end{equation*}
and
\begin{equation*}
\spec(\widehat A_0 + \widehat V\widehat X)\subset (-\infty,0].
\end{equation*}
Since the operator $ A_0 + V Y$ is diagonal with respect to
the decomposition $\fH_0 = \widehat{\fH}_0 \oplus \fK_0$,
\begin{equation*}
A_0 + V Y = (\widehat{A}_0 + \widehat{V} \widehat{X}) \oplus
(-VS|_{\fK_0}),
\end{equation*}
one infers that
\begin{equation}\label{sss}
\spec(A_0 + VY)\subset (-\infty,0].
\end{equation}

Combining \eqref{kkk}, \eqref{mmm}, and \eqref{sss} proves
that the operator $Y$ coincides with $X$ using the uniqueness
result of Theorem \ref{Clubshilfe}. Thus,  \eqref{vklu}
follows from \eqref{tak}.

Combining \eqref{vklu} and \eqref{kula} proves \eqref{ininin}.

The remaining statement of the theorem follows from the the
definition of the operator $Y$ and the fact that $X=Y$.
\end{proof}

By Theorem \ref{karkar} the uniqueness criterion (Theorem \ref{4.11a})
admits the following equivalent purely geometric formulation.

\begin{theorem}\label{4.11}
Assume Hypothesis \ref{diag} and let $X$ be the solution to the Riccati
equation
\begin{equation*}
A_1 X - X A_0 - X V X + V^\ast =0
\end{equation*}
referred to in Theorem \ref{Davis:Kahan}. Let $\fK_o$ and
$\fK_1$ be the subspaces given by \eqref{K0} and \eqref{K1},
respectively. Then  $X$ is the unique contractive solution if
and only if
\begin{itemize}
\item[(i)] either $\Ker(A_0-\lambda)\cap\Ker V^\ast$ or $\Ker(A_1-\lambda)\cap\Ker
V$ are trivial
\end{itemize}
and
\begin{itemize}
\item[(ii)] either $\fK_0$ or
$\fK_1$ (and hence both) are trivial.
\end{itemize}
The solution $X$ is strictly contractive.
\end{theorem}

\appendix
\section{Two Subspaces}\label{sec:2sub}

Here we collect some facts about two closed subspaces of a
separable Hilbert space which are used in the body of the paper.
Their comprehensive presentation with proofs as well as some
further results and the history of the problem can be found in
 \cite{Kostrykin:Makarov:Motovilov:2}.

Let $(P,Q)$ be an ordered pair of orthogonal projections in
the separable Hilbert space $\fH$. Denote
\begin{align*}
\fM_{pq}& :=\left\{f\in\fH\big| Pf=pf,\ Qf=qf \right\},\quad
p,q=0,1,
\\
\fM'_0& :=\Ran P\ \ominus\ (\fM_{10}\ \oplus\ \fM_{11}),
\\
\fM'_1& :=\Ran P^\perp\ \ominus\ (\fM_{00}\ \oplus\ \fM_{01}),
\\
\fM'& :=\fM'_0\oplus\fM'_1,
\\
P'& :=P|_{\fM'},
\\
Q'& :=Q|_{\fM'}.
\end{align*}
The space $\fH$ admits the canonical orthogonal decomposition
\begin{equation}\label{canon:decomposition}
\fH=\fM_{00}\ \oplus\ \fM_{01}\ \oplus\ \fM_{10}\ \oplus\
\fM_{11}\ \oplus\ \fM'.
\end{equation}

The following theorem provides a criterion for the subspace $\Ran Q$ to be
a graph subspace associated with the subspace $\Ran P$.
\begin{theorem}\label{3.2}
Let $P$ and $Q$ be orthogonal projections in a Hilbert space
$\fH$. The subspace $\Ran Q$ is a graph subspace $\cG(\Ran
P,X)$ associated with some closed densely defined (possibly
unbounded) operator $X: \Ran P\rightarrow \Ran P^\perp$ with
$\Dom(X)\subset\Ran P$ if and only if the subspaces
$\fM_{01}(P,Q)$ and $\fM_{10}(P,Q)$ in the canonical
decomposition \eqref{canon:decomposition} of the Hilbert space
$\fH$ are trivial, i.e.,
\begin{equation}\label{triv}
\fM_{01}(P,Q) = \fM_{10}(P,Q) = \{0\}.
\end{equation}
For a given orthogonal projection $P$ the correspondence
between the closed subspaces $\Ran Q$ satisfying \eqref{triv}
and closed densely defined operators $X:\Ran\to\Ran P^\perp$
is one-to-one.
\end{theorem}

The subspaces $\fM_{11}$ and $\fM_{00}$ have a simple
description in terms of the operator $X$: $\fM_{11}=\Ker X$
and $\fM_{00}=\Ker X^\ast$.

Note that $\fM_{01}(P,Q) = \fM_{10}(P,Q) = \{0\}$ if
$\|P-Q\|<1$. Moreover, Theorem \ref{3.2} has the following
corollary.

\begin{theorem}\label{frpr}
Let $P$ and $Q$ be orthogonal projections in a Hilbert space
$\fH$. Then the inequality $\|P-Q\|<1$ holds true if and only
if $\Ran Q$ is a graph subspace associated with the subspace
$\Ran P$ and some bounded operator $X\in\cB(\Ran P,\Ran
P^\perp)$, that is, $\Ran Q=\cG(\Ran P, X)$. In this case
\begin{equation}
\label{glgl} \|X\| = \frac{\|P-Q\|}{\sqrt{1-\|P-Q\|^2}}
\end{equation}
and
\begin{equation}
\label{estimate} \|P-Q\|=\frac{\|X\|}{\sqrt{1+\|X\|^2}}.
\end{equation}
\end{theorem}


\end{document}